\documentclass [twoside,10pt]{article}
\usepackage{graphicx}
\usepackage{amssymb}
\usepackage{amsthm}
\newtheorem{thm}{Theorem}
\newtheorem{cor}{Corollary}

\newtheorem{cla}{Claim}
\newcommand{\diam}{{\rm diam}}
\newcommand{\width}{{\rm width}}
\newcommand{\dist}{{\rm dist}}
\newcommand{\sh}{{\rm sinh \,}}
\newcommand{\ch}{{\rm cosh \,}}

\newcommand{\arch}{{\rm arccosh \,}}

\title{\bf Reduced polygons in the hyperbolic plane}

\date{}

\begin{document}

\maketitle

\thispagestyle{empty}

{ }
\vskip-1cm

\centerline{\author{MAREK LASSAK}}

\pagestyle{myheadings} \markboth{\centerline {Marek Lassak}}{\centerline {Reduced hyperbolic polygons}}

\baselineskip 12pt 

\maketitle
\vskip 0.5cm

\noindent 
{\bf Abstract}.
For a hyperplane $H$ supporting a convex body $C$ in the hyperbolic space $\mathbb{H}^d$ we define the width of $C$ determined by $H$ as the distance between $H$ and a most distant ultraparallel hyperplane supporting $C$. 
The minimum width of $C$ over all supporting $H$ is called the thickness $\Delta (C)$ of $C$.
A convex body $R \subset \mathbb{H}^d$ is said to be reduced if $\Delta (Z) < \Delta (R)$ for every convex body $Z$ properly contained in $R$.
We describe a class of reduced polygons in $\mathbb{H}^2$ and present some properties of them.
In particular, we estimate their diameters in terms of their thicknesses.

\baselineskip 16.7pt



\medskip

\date{}

\maketitle

\section{Introduction}

Let $H$ be a hyperplane supporting a convex body $C$ in the hyperbolic space $\mathbb{H}^d$. 
After \cite{[L23]} w define the {\bf width of $C$ determined by $H$} as the distance between $H$ and any farthest ultraparallel hyperplane supporting $C$. 
By compactness arguments there exists at least one such a most distant one, sometimes there are a finite or even infinitely many of them. 
The symbol ${\rm width}_H (C)$ denotes this width of $C$ determined by $H$.
Let us add that there also are different notions of width in $\mathbb{H}^d$.
For instance, these by Santal\'o, \cite{[Sa]}, Leichtweiss \cite{[Le]} and Horv\'ath \cite{[Ho]} who also presents a survey of earlier notions of width in $\mathbb{H}^d$.
Two more such surveys are given by B\"or\"oczky and Sagmeister in~\cite{[BS]} and by B\"or\"oczky, Cs\'epai and Sagmeister in
\cite{[BCS]}.

By {\bf the thickness $\Delta (C)$ of} a convex body $C \subset \mathbb{H}^d$ we mean the infimum of ${\rm width}_H (C)$ over all hyperplanes $H$ supporting $C$ (see \cite{[L23]}).
By compactness arguments, this infimum is realized, so $\Delta (C)$ is the minimum of the numbers ${\rm width}_H (C)$.

The above notions of width and thickness are analogous to the classic notions of width and thickness of a convex body in the Euclidean space $\mathbb{E}^d$ and in the spherical space
$\mathbb{S}^d$ (for the last one see \cite{[L15A]} and the survey article \cite{[L22]}) .

We say that $e$ is an {\bf extreme} point of a convex body $C \subset \mathbb{H}^d$ if the set $C \setminus \{ e \}$ is convex.

A convex body $R \subset \mathbb{H}^d$ is said to be {\bf reduced} if $\Delta (Z) < \Delta (R)$ for every convex body $Z$ properly contained in $R$.
This notion is analogous to the notion of a reduced body in $\mathbb{E}^d$ introduced by Heil \cite{[He]} (considered later also in the spherical space in  a number of papers, e.g., see \cite {[L15A]}--\cite{[L22]}).

In Section 2 we define a class of reduced polygons analogous to the classes of all reduced polygons in $\mathbb{E}^2$ and $\mathbb{S}^2$. 
It appears, that this analogy is not full.
Just this time there are also additional reduced polygons which do not fit to the patterns from $\mathbb{E}^2$ and $\mathbb{S}^2$.
Consequently, we call the polygons from the mentioned class ordinary reduced polygons. 
In Section 3 we present a number of properties concerning the width and thickness of ordinary reduced polygons.
Moreover, in Section 4 we estimate the diameter of the ordinary reduced polygons in terms of its thickness.
Section 5 is devoted to some questions on the diameter, perimeter, circumradius and inradius of ordinary reduced polygons.

For the  convenience of the reader, let us recall Proposition 1 (as Claim 1) and Theorem 1 from \cite{[L23]} (as Claim 2) which will be applied in this paper.

\begin{cla}
Let $C \subset \mathbb{H}^2$ be a convex body and $H$ be any supporting hyperplane of $C$.
Then $\width_H (C)$ equals to the maximum distance between $H$ and a point of $C$.
\end{cla}

\begin{cla}
For every convex body $C \subset \mathbb{H}^d$ we have 
\vskip -0.4cm

$$\max \{ {\rm width}_H (C); H {\rm \ is \ a \ supporting \ hyperplane \ of} \ C \} = {\rm  diam}(C).$$
\end{cla}

We interpret our considerations in the hyperboloid model of $\mathbb{H}^d$, so in the model on the upper sheet $x_{d+1}= \sqrt{x_1^2 + \dots + x_d^2 +1}$ of the two-sheeted hyperboloid.
This approach enables a proper similitude of the achieved results with the analogous facts in the Euclidean space and the spherical space. 
This model is considered by a number of mathematicians, for instance by Reynolds \cite{[Re]}.
Our figures present the orthogonal look to the sheet from the above.

For two points $p, q$, by $pq$ we denote the segment jointing them and by $|pq|$ its length.

The convex hull $V$ of $k \geq 3$ points in $\mathbb{H}^2$ such that each of them does not belong to the convex hull of the remaining points is called a {\it convex $k$-gon} in $\mathbb{H}^d$. 
These points are called {\bf vertices} of $V$. 
We write $V= v_1v_2\dots v_k$ provided $v_1, v_2, \dots , v_k$ are successive vertices of $V$, when we go around $V$ on the boundary of $V$ according to the positive orientation. 
When we take $k \geq 3$ points $v_1, \dots , v_k$ in a hyperbolic circle of $\mathbb{H}^2$ such that $|v_1v_2| = ... = |v_{k-1}v_k| = |v_kv_1|$, the convex hull of them is called a {\bf hyperbolic regular $k$-gon}.

\section{On some reduced odd-gons}

In a convex $n$-gon $v_1 \dots v_n$ we consider indices modulo $n$. 

The following theorem is an analog of the ``if" parts of Theorem 7 of \cite{[L90]} and Theorem 3.2 of \cite{[L15C]}. 

\begin{thm} \label{ShapeV} 
Let $V = v_1 \dots v_n$ be a convex odd-gon.
Assume that for every  $i \in \{1, \dots , n\}$ the projection of $v_i$ on the line $L_i$ containing the side $S_i = \ v_{i+(n-1)/2} v_{i+(n+1)/2}$ is in the relative interior of this side
and the distance between $v_i$ and $L_i$ is $\Delta (V)$. 
Then $V$ is reduced.
\end{thm}

\begin{proof}
Consider an odd-gon $V = v_1v_2\dots v_n$ such that the projection of every $v_i$ onto the straight line $L_i$ containing the side $S_i = \ v_{i+ (n-1)/2}v_{i+ (n+)/2}$ is in the relative interior of this side and that all the distances $\dist(L_i, v_i)$ are equal.
Since the projections of vertices $v_{i+ (n-1)/2}, v_i, v_{i+ (n+1)/2}$, where $i \in \{1, \dots ,n \}$, onto the lines containing the ``opposite" sides $S_{i+ (n-1)/2}, S_i, S_{i+ (n+1)/2}$, respectively, are in the relative interiors of these sides, the vertex $v_i$ is the only vertex in the distance $\Delta(V)$ from the straight line $L_i$ containing $S_i$ (so the remaining are closer, i.e, strictly between $L_i$ and the equidistance curve to $L_i$ ``supporting" $V$).
This is true for every $i \in \{1,\dots, n\}$.
Consequently, the fact that all $\dist(L_i, v_i)$ are equal imply the following. 
If we take an arbitrary convex body $Z \subset V$ different from $V$, then $Z$ does not contain a vertex of $V$. 
Consequently, $\Delta (Z) < \Delta (V)$. 
Thus $V$ is reduced.
\end{proof}

Let us call the polygons described in this theorem as {\bf ordinary reduced polygons} in $\mathbb{H}^2$.
The reason of this name is that the constructions of them are analogous to all reduced polygons in $\mathbb{E}^2$ and $\mathbb{S}^2$ (again see Theorem 7 of \cite{[L90]} and Theorem 3.2 of \cite{[L15C]}),
and that, surprisingly, there exist some reduced polygons that are not ordinary reduced polygons.
Namely, the author recently learned and checked that some hyperbolic rhombi (i.e., convex hulls of two perpendicular segments intersecting each other at midpoints) are reduced. 
The proof of this statement is a part of an ongoing project on hyperbolic reducedness by K. Jr. B\"or\"oczky, A. Freyer and \'A. Sagmeister.
By the way, analogously some crosspolytopes in $\mathbb{H}^d$ (i.e., convex hulls of $d$ perpendicular segments intersecting each other at midpoints) are reduced polytopes. 
This fact is interesting since still an open question is if there are reduced polytopes in $\mathbb{E}^d$ for $d >3$ (see p. 374 of \cite {[L90]} and Problem 5 of \cite{[LM]}). 
They exist for $d=3$ as it follows from the paper \cite{[GJPW]} by Gonzalez Merino, Jahn, Polyanskii and Wachsmuth.

We see some ordinary reduced pentagon and heptagon in Fig. 1 and Fig. 2, respectively.

\vskip0.2cm

\begin{figure}[htbp]        
\hskip2.75cm
\includegraphics[width=8.4cm, height=8.24cm]{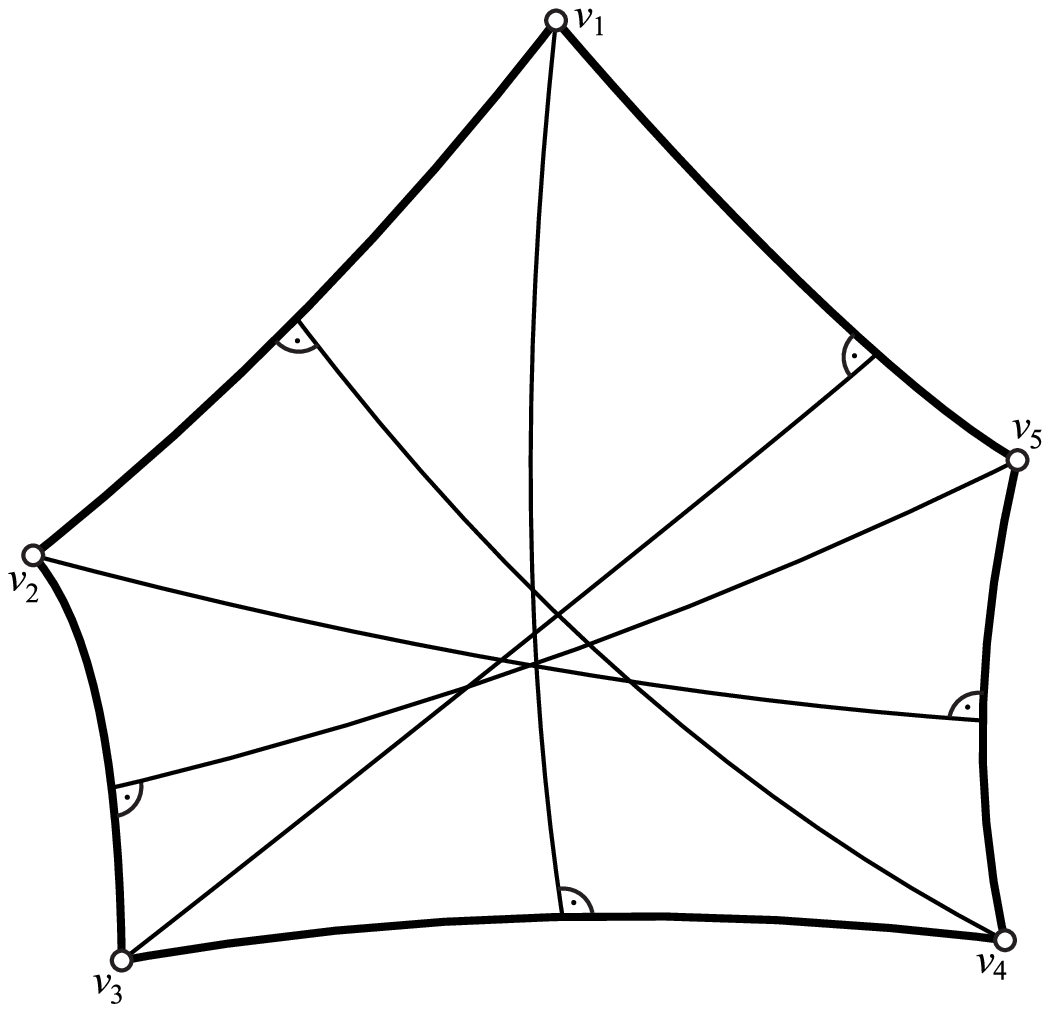}

\vskip0.3cm 
\centerline
{{\bf Fig. 1.} An ordinary reduced pentagon} 
\end{figure}

\begin{cor}  
Every regular hyperbolic odd-gon is reduced. 
\end{cor}

\begin{cor} 
The only reduced hyperbolic triangles are the regular ones. 
\end{cor}

The first corollary is obvious.
In order to show the second one imagine that a non-regular reduced triangle $T$ exists.
Its three altitudes are not equal.
Take any shortest. 
Its length is equal to $\Delta(T)$.
Denote by $z$ the vertex at a longer altitude.
By a straight line, we cut off from $T$ a sufficiently small piece containing $z$.
We obtain a convex quadrangle of thickness $\Delta(T)$ in contradiction to the assumption that $T$ is reduced.

\begin{cor} 
Assume that a straight line $L$ supports an ordinary reduced polygon $V$.
We have $\width_L (V) = \Delta (V)$ if and only if $L$ contains a side of $V$.
\end{cor}

\section{Some properties of ordinary reduced polygons}

Here is a lemma for $\mathbb{H}^2$ (true also in $\mathbb{E}^2$ and $\mathbb{S}^2$) which is needed in the proof of part (iii) of the forthcoming theorem. 

\vskip0.25cm
\noindent
{\bf Lemma.}
{\it If $U \subset W$ are different convex polygons, then the perimeter of $U$ is smaller than this of $V$.}

\vskip0.25cm
For the proof observe that we can get $V$ from $W$ by a finite number of successive cuttings of non-degenerate triangles.
By the triangle inequality, each time the perimeter decreases.
This implies our thesis.

In the following theorem we use the notation from Theorem 1.
Moreover, by $p_i$ denote the projection of $v_i$ onto the line $L_i$ containing the segment $v_{i+(n-1)/2} v_{i+(n+1)/2}|$.

\begin{thm} \label{properties}
Let $V$ be any ordinary reduced odd-gon $v_1 \dots v_n$.
Then

{\rm (i)} we have $|v_ip_{i+(n+1)/2}| = |p_iv_{i+(n+1)/2}|$ for $i=1, \dots , n$,

{\rm (ii)} for every $i \in \{1, \dots , n\}$ the segment $v_ip_i$ halves the perimeter of $V$,

{\rm (iii)} if $V$ is different from a triangle, then we have $\beta_i < \alpha_i$ for $i= 1, \dots ,n$, where $\alpha_i = \angle v_{i+1}v_ip_i$ and $\beta_i = \angle p_iv_ip_{i+(n+1)/2}$, and for $V$ as a triangle we have $\beta_i = \alpha_i$. 
\end{thm}

\begin{proof}
(i) For every $i \in \{1, \dots , n\}$ take the line $L_i$ containing the side $v_{i+(n-1)/2}v_{i+(n+1)/2}$ (see Fig. 2). 
By the definition of the ordinary reduced polygons, the projection $p_i$ of $v_i$ onto $L_i$ belongs to the relative interior of this side and the distance of $v_i$ from $L_i$, so $|v_ip_i|$ is $\Delta(V)$.

\begin{figure}[htbp]        
\hskip0.05cm
\includegraphics[width=14.345cm, height=6.555cm]{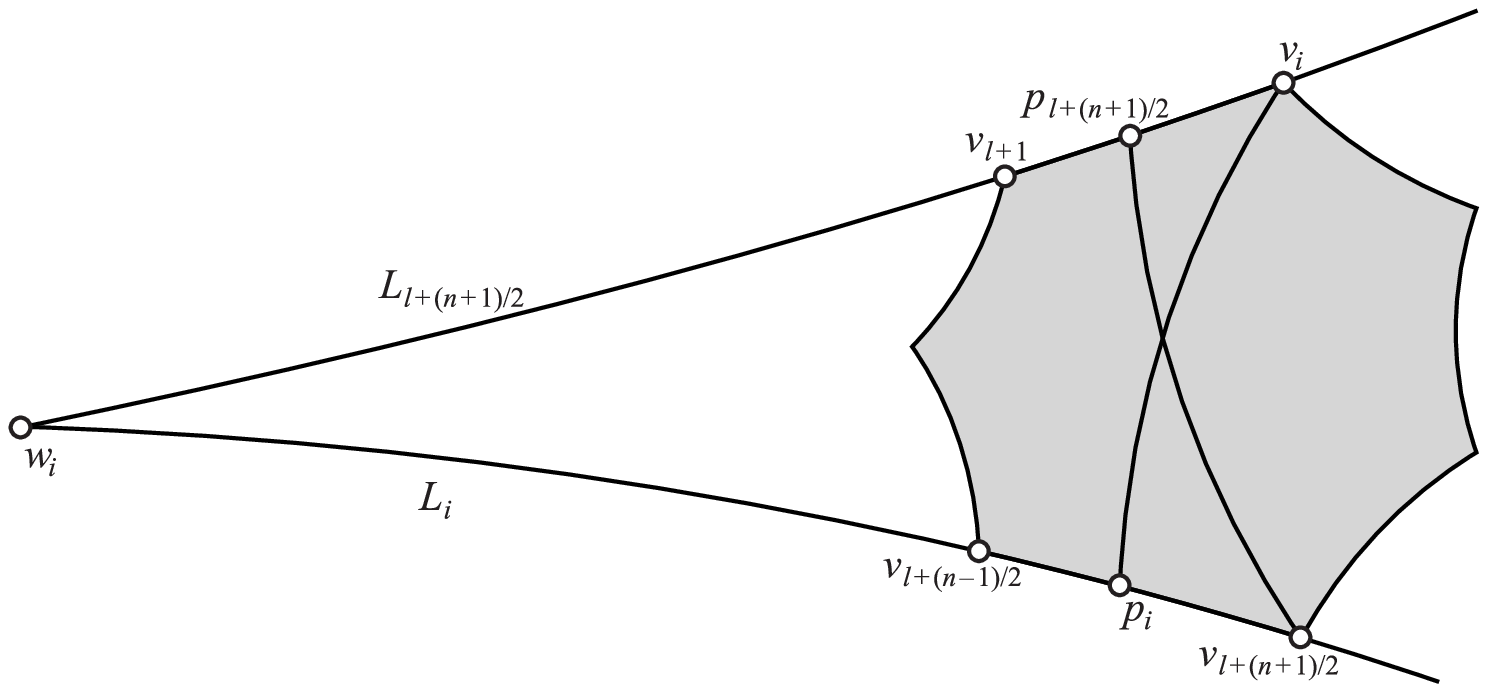}

\vskip0.4cm 
\centerline
{{\bf Fig. 2.} Illustration to the proof of Theorem \ref{properties}} 
\end{figure}

Take the line $L_{i+(n+1)/2}$ containing the side $v_iv_{i+1}$.
We know that the projection $p_{i+(n+1)/2}$ of $v_{i+(n+1)/2}$ onto $L_{i+(n+1)/2}$ belongs to the relative interior of this side.
We also know that the distance of $v_{i+(n+1)/2}$ from $L_{i+(n+1)/2}$, so $|v_{i+(n+1)/2}p_{i+(n+1)/2}|$ is $\Delta(V)$.

By $w_i$ denote the intersection point of $L_i$ and $L_{i+(n+1)/2}$ if they intersect (\'Ad\'am Sagmeister payed my attention 
that this is not always the case, thanks!)
Since $|v_{i+(n+1)/2} p_{i+(n+1)/2}| = \Delta(V) = |v_ip_i|$, the rectangular triangles $w_ip_{i+(n+1)/2}v_{i+(n+1)/2}$ and $w_ip_iv_i$ are congruent (they are reflected copies of each other about the angular bisector of the triangle $w_iv_iv_{i +(n+1)/2}$ at $w_i$).
Hence $|w_iv_{i+(n+1)/2}| = |w_iv_i|$ and $|w_ip_{i+(n+1)/2}| = |w_ip_i|$.
Consequently, $|v_ip_{i+(n+1)/2}| = |p_iv_{i+(n+1)/2}|$.

If $L_i$ and $ L_{i+(n+1)/2}$ do not intersect, provide a straight line that intersects both $L_i$ and $L_{i+(n+1)/2}$ at points $a_i \in L_i$ and $b_i \in  L_{i+(n+1)/2}$ with equal angles $p_ia_ib_i$ and $a_ib_ip_{i+(n+1)/2}$ to both these lines. 
Moreover let $p_i$ belong to the relative interior of $a_iv_{i+(n+1)/2}$ and let $p_{i+(n+1)/2}$ belong to the relative interior of $b_iv_i$. 
Observe that the quadrangles $v_ip_ia_ib_i$ and $v_{i+(n+1)/2}p_{i +(n+1)/2}b_ia_i$ are congruent (still they have equal corresponding angles, the common side $a_ib_i$ and $|v_ip_i| = |v_{i+(n+1)/2}p_{i+(n+1)/2}|$).
Hence (i) holds true.

(ii) Applying $n$ times (i) we conclude that the sum of the lengths of the boundary segments of $V$ from $v_i$ to $p_i$ (moving according to the positive orientation) is equal to the sum of the boundary segments from $p_i$ to $v_i$ (again moving according to the positive orientation).
Therefore, we obtain the required thesis.

(iii) First assume that $L_i$ and $L_{i+(n+1)/2}$ do intersect.
By (i) and (ii) the lengths of the pieces $B_i$ of the boundary of $V$ from $p_{i + (n+1)/2}$ to $p_i$ and $B'_i$ from $v_{i + (n+1)/2}$ to $v_i$ are equal.
Since $B_i$ is in the triangle $T_i= p_{i + (n+1)/2}w_ip_i$, then by Lemma for the triangle $T_i$ as $W$ and the polygon $p_{i + (n+1)/2}v_{i+1}\dots v_{i + (n-1)/2}p_i$ as $U$ we get that the length of $B_i$ is below $|p_{i + (n+1)/2}w_i| + |w_ip_i|$ (recall that $w_i$ is defined in (i) and see Fig. 2).
Analogously, the length of the boundary of $V$ from $v_{i + (n+1)/2}$ to $v_i$ is over $|v_iv_{i + (n+1)/2}|$.
Hence $|v_{i + (n+1)/2}v_i| < |p_{i + (n+1)/2}w_i| + |w_ip_i|$.
From (i) we obtain $|p_iv_{i + (n+1)/2}| + |v_iv_{i + (n+1)/2|} < |p_iw_i| + |w_iv_i|$.
Since the triangles $v_ip_iv_{i + (n+1)/2}$ and $v_iw_ip_i$ have the common side $v_ip_i$ and right angles at $p_i$, 
we get $|p_ip_{i + (n+1)/2}| < |w_ip_i|$, implying $\beta_i < \alpha_i$.
Clearly, if $V$ is a triangle, then by Corollary 2 it is a regular triangle and thus $\beta_i = \alpha_i$.

If $L_i$ and $L_{i+(n+1)/2}$ do not intersect, we again deal with the quadrangles, as in the last paragraph of the proof of (i).
Applying the observation that the length of $B_i$ is below $|p_{i+(n+1)/2}b_i| + |b_ia_i| + |a_ip_i|$, we analogously obtain that $\beta_i < \alpha_i$.
\end{proof}

\section{On the diameter of ordinary reduced polygons}

\noindent
{\bf Proposition.}
{\it The diameter of any ordinary reduced $n$-gon is realized for only for some pairs of vertices whose indices (modulo $n$) differ by $(n-1)/2$ or by $(n+1)/2$.} 

\begin{proof}
Take an ordinary reduced $n$-gon $V =v_1 \dots v_n$.
Clearly its diameter is realized for a pair of vertices.
We have to show that the diameter is $|v_iv_{i+(n-1)/2}|$ or $|v_iv_{i+(n+1)/2}|$ for an $i \in \{1, \dots , n\}$.
 
Claim 1 says that $\diam(V)$ equals to the maximum $\width_L (V)$ over all straight lines $L$ supporting $V$ and Claim 2 states that $\width_L (V)$ is equal to the maximum distance of a point, so a vertex, of $V$ from $L$. 
Hence $\diam(V)$ equals to the maximum distance from $L$, over all supporting straight lines $L$, to a farthest vertex of $V$.

Let us continuously change the position of the supporting line $L$ of $V$ according to the positive orientation from the position in which it contains a side $S_i$ up to the position when it contains the next side $S_{i+1}$, all the time looking for the farthest points of $V$ from $L$. 
During this changing, the farthest points of $V$ from $L$ may be only $v_i$ or $v_{i+1}$. 
More precisely, having in mind the definition of an ordinary reduced polygon, we see that at the beginning positions of $L$ the most distant vertex from each such $L$ is only $v_i$, then for exactly one position of $L$ both $v_i$ and $v_{i+1}$ are most distant from $L$, and finally up to the moment when $L$ contains $S_{i+1}$ only $v_{i+1}$ is the most distant vertex to $L$.
Consequently, we see that the most distant vertices from $v_i$ may be only $v_i$ or $v_{i+1}$.

We conclude that the diameter of $V$ is $|v_iv_{i+(n-1)/2}|$ or $|v_iv_{i+(n+1)/2}|$ for an $i \in \{1, \dots , n\}$.
\end{proof}

The following theorem is analogous to the second part of Theorem 9 from \cite{[L90]} and Theorem 3.1 from \cite{[L15C]} (for completeness let us also mention the paper \cite{[LC]} by Liu and Chang).
The proof is analogous, but we provide it here, since we apply some facts established in the present paper and some formulas from the hyperbolic plane.

\begin{thm} \label{diameter}
For every ordinary reduced polygon $V \subset \mathbb{H}^2$ we have

$$\diam (V) < \arch \Big( \sqrt {1 + {1 \over 3} \sh \Delta (V)} \cdot \ch \Delta (V) \Big).$$
\end{thm}

\begin{proof}
By Proposition there is an $i \in \{1, \dots , n\}$ such that $\diam (V)$ equals to $|v_iv_{i+(n+1)/2}|$ or   $|v_iv_{i+(n-1)/2}|$.
Consider the first possibility (in the second one, the consideration is analogous).

Put $r_i = |p_iv_{i + (n+1)/2}|$, $s_i = |v_iv_{i + (n+1)/2}|$ and $\gamma_i = \angle p_iv_{i + (n+1)/2}v_i$.  

By $o_i$ denote the intersection point of $v_ip_i$ and $v_{i + (n+1)/2}p_{i + (n+1)/2}$.
From the observation in the brackets in the last paragraph of the proof of (i) in Theorem \ref{properties} we conclude that $|v_io_i| = |v_{i + (n+1)/2}o_i|$.
Hence the triangle $v_io_iv_{i + (n+1)/2}$ is isosceles.
Thus $\angle o_iv_{i + (n+1)/2}v_i = \angle o_iv_{i + (n+1)/2}v_i$.
So $\gamma_i = \alpha_i + \beta_i$.
This and (iii) in Theorem \ref{properties} imply $\gamma_i \geq 2\beta_i$.
From $|v_ip_i| = \Delta (V)$ and by the hyperbolic law of sines for the right triangle $v_ip_iv_{i + (n+1)/2}v_i$ we obtain ${\sh r_i \over \sin \beta_i} = {\sh \Delta (V) \over \sin \gamma_i}$.
Thus by $\gamma_i \geq 2\beta_i$ we get ${\sh r_i \over \sin \beta_i} \leq {\sh \Delta (V) \over \sin 2 \beta_i}$.
This and $\sin 2\beta_i = 2\sin \beta_i \cos\beta_i$ imply 

\vskip-0.1cm
$$\sh r_i \leq {\sh \Delta (V) \over 2\cos \beta_i}.$$ 

\vskip0.15cm
Let us show that $\beta_i < {\pi \over 6}$.
Imagine the opposite case when $\beta_i \geq {\pi \over 6}$.
Then from $\beta_i \leq \alpha_i$ (shown in (iii) of Theorem \ref{properties}) and $\gamma_i = \alpha_i + \beta_i$ it follows $\gamma_i \geq 2\beta_i > \frac{\pi}{3}$ and hence $\gamma_i \geq \frac{\pi}{3}$, which leads to the conclusion that the sum of angles in the triangle $v_ip_iv_{i +(n+1)/2}$ is at least $\pi$. 
This in the hyperbolic triangle is impossible.
So $\beta_i < {\pi \over 6}$.

From $\beta_i < {\pi \over 6}$ we get $\cos \beta_i > \frac{\sqrt 3}{2}$.
Thus by the inequality just before the preceding paragraph we obtain $\sh r_i \ < \frac{\sh \Delta (V)}{\sqrt 3}$.
Hence $\sh^{\hskip-0.06cm 2} r_i < \frac{\sh^{\hskip-0.06cm 2} \Delta (V)}{3}$.
After we apply $\ch^{\hskip-0.06cm 2} r_i - \sh^{\hskip-0.06cm 2} r_i = 1$, i.e., $\sh^{\hskip-0.06cm 2} r_i = \ch^{\hskip-0.06cm 2} r_i -1$, we obtain $\ch^{\hskip-0.06cm 2} r_i < 1 + \frac{1}{3} \sh^{\hskip-0.06cm 2} \Delta (V)$.
So $\ch r_i < \sqrt{1+\frac{1}{3} \sh^{\hskip-0.06cm 2} \Delta (V)}$.

Take into account the equality $\ch s_i = \ch r_i \ch \Delta (V)$, which follows from the hyperbolic Pythagorean theorem for the triangle $v_ip_iv_{i+ (n+1)/2}$.
In other words, $\ch r_i = \frac{\ch s_i}{\ch \Delta (V)}$.

By the preceding two paragraphs we get the following inequality implying the next one

$$\frac{\ch s_i}{\ch \Delta (V)} < \sqrt{1 + \frac{1}{3} \sh^{\hskip-0.06cm 2} \Delta (V)},$$ 

$$s_i < \arch \Big(\sqrt{1+ \frac{1}{3} \sh^{\hskip-0.06cm 2} \Delta (V)} \cdot \ch \Delta (V)\Big).$$ 

Consequently, by Proposition we obtain the thesis.
\end{proof}

\section{A few questions}

The author expects that $1 < \frac{\diam (V)}{\Delta (V)} < 2$ for every ordinary reduced polygon $V \subset \mathbb{H}^2$ with $2$ as the limit of this quotient for the regular triangle when $\Delta (V)$ tends to $\infty$.
On the other hand, $1$ is the limit of this quotient for the limit of the sequence of regular odd-gons of a fixed thickness.
This follows from $\ch r_i = \frac{\ch s_i}{\ch \Delta (V)}$ established in the proof of Theorem \ref{diameter} and from Proposition. 
Still $s_i$ is the diameter of the regular $n$-gon with $n$ odd and $r_i \to 0$.

We conjecture that the diameter of any non-regular ordinary reduced $n$-gon in $\mathbb{H}^2$ is greater than the diameter of the regular $n$-gon of the same thickness, analogously as in $\mathbb{E}^2$ and on $\mathbb{S}^2$.
For $\mathbb{E}^2$ this is formulated in Corollary 6 of \cite {[L90]} and for $\mathbb{S}^2$ this is proved in Theorem 4.1 from the paper \cite{[CHJ]} by Chen, Hou and Jin.

We may also consider the question about minimizing the diameter for general ordinary convex polygons of a fixed thickness in $\mathbb{H}^2$, equivalent to maximizing the thickness for convex polygons of a fixed diameter.
For this compare the paper \cite{[BF]} by Bezdek and Fodor on maximizing the thickness (called ``width" there) for the unit diameter convex polygons in $\mathbb{E}^2$.
Also see the related paper \cite{[AHM]} by Audet, Hansen and Messine.

Is the perimeter of any ordinary non-regular reduced $n$-gon in $\mathbb{H}^2$ larger than the perimeter of the regular $n$-gon?
The same question for the ``area" in place of the ``perimeter".

Recall that every reduced polygon $R \subset E^2$ is contained in a disk of radius $\frac{2}{3} \Delta(R)$ as shown in  Proposition from \cite{[L03]}.
This generates the following problem for $\mathbb{H}^2$.
What is the smallest radius of a disk which contains every ordinary reduced polygon of a given thickness in $\mathbb{H}^2$?
Let us add that an analogous question for a reduced polygon on $\mathbb{S}^2$ is solved in \cite{[LC+]} by Liu and Chang.
By the way, why do not consider the dual problem on the inradius of an ordinary reduced polygon in $\mathbb{H}^2$.
The same question for the reduced polygons in $\mathbb{E}^2$ and $\mathbb{S}^2$.

\baselineskip 12pt

\vskip0.1cm
\noindent
MAREK LASSAK

\noindent
University of Science and Technology

\noindent
85-789 Bydgoszcz, Poland

\noindent
e-mail: lassak@pbs.edu.pl

\end{document}